\documentclass[11pt]{article}
\title{On Characteristic Polynomials of the Family of Cobweb Posets}
\author{Ewa Krot-Sieniawska \\
\\Institute of Computer Science, Bia{\l}ystok University\\
PL-15-887 Bia{\l}ystok, ul.Sosnowa 64, POLAND\\
e-mail: ewakrot@wp.pl, ewakrot@ii.uwb.edu.pl}
\date{}

\usepackage{amsmath,amsthm}
\usepackage{graphicx}
\chardef\bslash=`\\ 
\newtheorem{ex}{Example}[section]
\newtheorem{defn}{Definition}[section]

\newtheorem{thm}{Theorem}[section]
\newtheorem{prop}{Proposition}[section]
\newtheorem{cor}{Corollary}[section]

\begin{document}
\maketitle
\begin{abstract}
\noindent This note is a response to one of problems posed by A.K.
Kwa\'sniewski in \cite{49}. Namely with $\{P_n\}_{n\geq 0}$ being the
sequence of finite cobweb subposets, the looked for explicit formulas for
corresponding sequence $\{\chi_n(t)\}_{n\geq 0}$ of $P_n$'s
characteristic polynomials are discovered and delivered here. The recurrence
relation defining  arbitrary family $\{\chi_n(t)\}_{n\geq 0}$ is also
derived.

\end{abstract}
AMS Classification numbers: 06A06,  06A07, 06A11,  11C08, 11B37\\

\noindent Presented at Gian-Carlo Rota Polish Seminar:
http://ii.uwb.edu.pl/akk/sem/sem rota.htm
\section{ Cobweb posets}
The   family of the so called cobweb posets $\Pi$ has been
invented by A.K.Kwa\'sniewski few years ago (for references  see:  \cite{44,46}).
These structures are such a generalization of the Fibonacci tree growth  that
allows joint combinatorial interpretation for all of them under
the admissibility condition (see \cite{49,49a}).

\noindent Let $\{F_n\}_{n\geq 0}$ be a natural numbers valued sequence with
$F_0=1$ (with $F_0=0$ being exceptional as in case of Fibonacci
numbers). Any sequence satisfying this property uniquely
designates cobweb poset  defined as follows.

\noindent For $s\in\bf{N}_0=\bf{N}\cup\{0\}$ let us to define levels of
 $\Pi$:
$$\Phi_{s}=\left\{\langle j,s \rangle ,\;\;1\leq j \leq F_{s}\right\},\;\;\;$$
(in case of $F_0=0$  level $\Phi_0$ corresponds to the empty root
$\{\emptyset\}$). )

\noindent Then

\begin{defn}
Corresponding cobweb poset is  an infinite partially ordered set
$\Pi=(V,\leq)$, where
 $$ V=\bigcup_{0\leq s}\Phi_s$$
 are the elements ( vertices) of $\Pi$ and the partial order relation $\leq$ on $V$ for
 $x=\langle s,t\rangle, y=\langle u,v\rangle $ being  elements of
cobweb poset $\Pi$ is defined by  formula
$$ ( x \leq_{P} y) \Longleftrightarrow
 [(t<v)\vee (t=v \wedge s=u)].$$
\end{defn}

\noindent Obviously any cobweb poset can be represented, via its Hasse
diagram, as infinite directed  graf  $\Pi=\left( V,E\right)$,
where  set $V$ of its vertices is defined as above and

$$E =\{\left(\langle j , p\rangle,\langle q ,(p+1) \rangle
\right)\}\;\cup\;\{\left(\langle 1 , 0\rangle ,\langle 1 ,1
\rangle \right)\},$$ \quad where $1 \leq j \leq {F_p}$ and $1\leq
q \leq {F_{(p+1)}}$ stays for  set of (directed) edges.

\noindent For example the Hasse diagram of Fibonacci cobweb poset designated
by the famous Fibonacci sequence is presented below.

\begin{center}
\includegraphics[width=100mm]{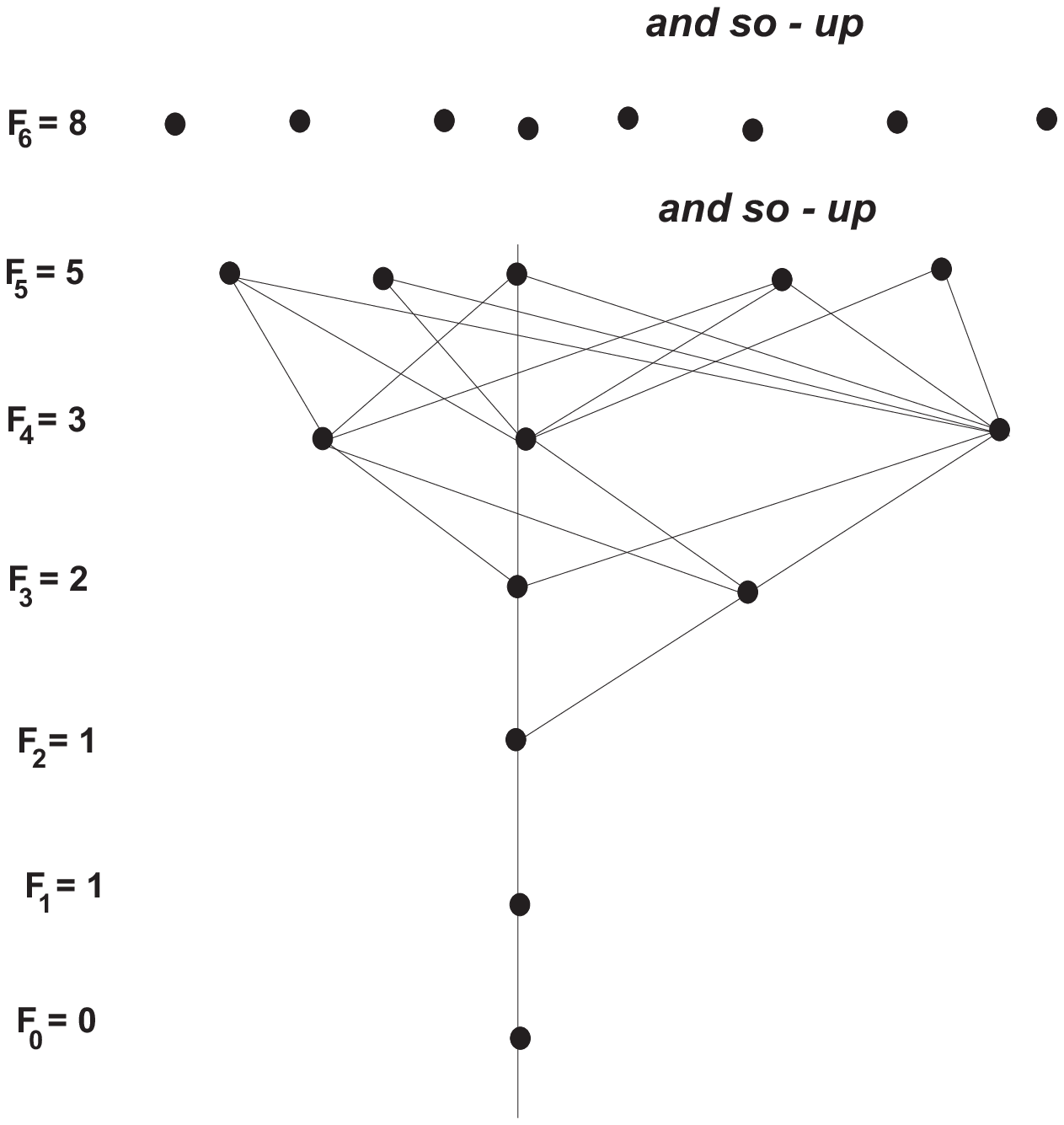}

\vspace{2mm}

\noindent {\small Fig.~1. The construction of the  Fibonacci
cobweb  poset}
\end{center}

\noindent The Kwasniewski cobweb posets under consideration represented by  graphs  are examples of
oderable directed acyclic graphs (oDAG)  which we start to call from now in brief:  KoDAGs. These
are  structures of universal importance for the whole of mathematics - in particular for
discrete "`mathemagics"' [http://ii.uwb.edu.pl/akk/ ]  and computer sciences in general
(quotation from \cite{49,49a} ):

\begin{quote}
For any given natural numbers valued sequence the graded (layered)
cobweb posets` DAGs  are equivalently representations of a chain
of binary relations. Every relation of the cobweb poset chain is
biunivocally represented by the uniquely designated
\textbf{complete} bipartite digraph-a digraph which is a
di-biclique  designated  by the very  given sequence. The cobweb
poset is then to be identified with a chain of di-bicliques i.e.
by definition - a chain of complete bipartite one direction
digraphs.   Any chain of relations is therefore obtainable from
the cobweb poset chainof complete relations  via
deleting  arcs (arrows) in di-bicliques.\\
Let us underline it again : \textit{any chain of relations is
obtainable from the cobweb poset chain of  complete relations via
deleting arcs in di-bicliques of the complete relations chain.}
For that to see note that any relation  $R_k$ as a subset of  $A_k
\times A_{k+1}$ is represented by a  one-direction bipartite
digraph  $D_k$.  A "complete relation"  $C_k$ by definition is
identified with its one direction di-biclique graph $d-B_k$.  Any
$R_k$ is a subset of  $C_k$. Correspondingly one direction digraph
$D_k$ is a subgraph of an one direction digraph of $d-B_k$.\\
The one direction digraph of  $d-B_k$ is called since now on
\textbf{the di-biclique }i.e. by definition - a complete bipartite
one direction digraph.   Another words: cobweb poset defining
di-bicliques are links of a complete relations' chain.
\end{quote}

\noindent According to the  definition above arbitrary cobweb
poset $\Pi=(V,\leq)$ is a graded poset ( ranked poset) and for $s\in\bf
 N_0$:
 $$x\in\Phi_s\;\; \longrightarrow\;\; r(x)=s,$$
 where $r:\Pi\rightarrow \bf N_0$ is a rank function on $\Pi$.

\noindent Let us then define Kwa\'sniewski finite cobweb
sub-posets as follows
 \begin{defn}
 Let $P_n=(V_n,\leq)$, $(n\geq 0)$, for ${\displaystyle
V_n=\bigcup_{0\leq s\leq n}}\Phi_s$ and  $\leq$ being the induced
partial order relation on $\Pi$.
\end{defn}

\noindent Its easy to see that $P_n$ is ranked poset with rank function $r$
as above. $P_n$ has a unique minimal element $0=\langle
1,0\rangle$ ( with $r(0)=0$). Moreover $\Pi$ and all $P_n$ s
satisfy the Jordan chain condition and the length of $P_n$ is
$l(P_n)=r(P_n)=n$ for $n\geq 0$.

\noindent For finite graded poset $P$ one can define (see \cite{29g})
Whitney numbers of the first and second kind $w_k(P)$ and $W_k(P)$
respectively as follows
$$w_{k}(P)=\sum_{x\in P,\,r(x)=k}\mu(0,x),$$
$$W_k(P)=\sum_{x\in P,\,r(x)=k}1=|\{x\in P\,:\,r(x)=k\}|,$$
where $\mu$ stays for M\"{o}bius  function of $P$ indispensable in
numerous inversion type formulas of countless applications (see
\cite{29g, 65, b, 73,74}).

\noindent Then the characteristic polynomial of $P$ \cite{65, b, 73,74} is
the polynomial
$$\chi_{P}(t)=\sum_{x\in P}\mu(0,x)t^{n-r(x)}=\sum_{k= 0}^nw_k(P)t^{n-k},$$
where $n=l(P)$.

\noindent Here next we answer the question posed by A.K.Kwa\'sniewski in the source
paper for the problem in question \cite{49}.
\begin{quote}
{\em Let $\{P_n\}_{n\geq 0}$ be the sequence of finite cobweb
subposets (...). What is the form and properties of
$\{P_n\}_{n\geq 0}$'s characteristic polynomials
$\{\rho_n(\lambda)\}_{n\geq 0}$? (...) What are recurrence
relations defining the  family $\{\rho_n(\lambda)\}_{n\geq 0}$?}
\end{quote}

\section{Whitney numbers of cobweb posets}
Obviously for arbitrary cobweb poset  $\Pi$ and for all its finite
subposets  $P_n$, ($n\geq 0$) one has:
\begin{equation}
W_k(\Pi)=F_k,\;\;\;\;\;k\geq 0,
\end{equation}
where $\{F_n\}_{n\geq 0}$ is a natural numbers valued sequence
uniquely designating $\Pi$.

\noindent Now let us consider the corresponding  numbers $w_k(\Pi)$.
The explicite formula for M\"{o}bius function of the Fibonacci
cobweb poset uniquely designated by  the  Fibonacci sequence  was
derived by the present author in \cite{35, 36}. It can be easy extend to the hole family of cobweb posets and their
finite subposets $P_n$, ($n\geq 0$),\cite{37a}. Moreover, by the use of notion of the standard reduced incidence algebra $R(\Pi)$,  (see \cite{37a}) one can  show, that for $x\in\Pi$ the value $\mu(0,x)$ depends on $r(x)$
only. So for $x$ as above we have:

\begin{equation}
\mu(0,x)=\mu(r(0),r(x))=\mu(0,r(x))=\mu(r(x)).
\end{equation}
Moreover
\begin{equation}
\mu(0,x)=\mu(r(x))=(-1)^{r(x)}\prod_{i=1}^{r(x)-1}(F_i-1).
\end{equation}

\noindent Then \begin{prop} For arbitrary cobweb poset $\Pi$ and for all its
finite subposets  $P_n$, ($n\geq 0$) corresponding Whitney numbers
of the first kind are given by the formulas:

\noindent for  $k>0$
\begin{equation}
\begin{array}{lll}
  w_k(\Pi) & = & {\displaystyle \sum_{\{x\in\Pi\,:\,r(x)=k\}}}\mu(0,x)=F_k\cdot\mu(0,x) \\
   &  &  \\
   & = & F_k\cdot(-1)^k \cdot{\displaystyle \prod_{i=1}^{k-1}}(F_i-1) \\
\end{array}
\end{equation}
and
\begin{equation}
w_0(\Pi)=1.
\end{equation}
\end{prop}

\section{The characteristic polynomials of finite cobweb posets}

\noindent The knowledge of Whitney numbers $w_k(P_n)$,  enables us to
construct the characteristic polynomials for all $P_n$, ( $n\geq
0$). Let us recall the formula defining $\chi_n(t)$:
$$\chi_{P_n}(t)=\sum_{x\in P_n}\mu(0,x)t^{n-r(x)}=\sum_{k= 0}^nw_k(P_n)t^{n-k}.$$
Using the above formulas one has

\begin{thm} The characteristic polynomials  $\chi_{P_n}(t)$, ($n\geq 0$) are given by the following explicit formula:
\begin{equation}
\chi_{P_n}(t)=\chi_n(t)=x^n+\sum_{k=1}^n (-1)^kF_k\cdot
\prod_{i=1}^{k-1}(F_i-1)x^{n-k}.
\end{equation}
\end{thm}
Moreover, as in case of Fibonacci cobweb poset, the following holds:

\begin{cor} Let $\{F_n\}_{n\geq}$ be the sequence designating the
cobweb poset $\Pi$ (and all corresponding sub-posets $P_n$). In
the case $F_1=1$ ( or equivalently $|\Phi_1|=1$) one has
\begin{equation}
\chi_{n}(t)=t^n-t^{n-1}
\end{equation}
for $n\geq 1$ and
\begin{equation}
\chi_0(t)=1.
\end{equation}
\end{cor}
\begin{cor}
Let $\{F_n\}_{n\geq}$ be the sequence designating the cobweb poset
$\Pi$ (and all corresponding sub-posets $P_n$). Then the sequence
$\{\chi_n(t)\}_{n\geq 0}$ of $\{P_n\}_{n\geq 0}$'s characteristic
polynomials is defined by the following recurrence relation
\begin{equation}
\chi_0(t)=1, \quad\quad \chi_1(t)=t-F_1
\end{equation}
\begin{equation}
\chi_n(t)=t\chi_{n-1}(t)+(-1)^nF_n(F_{n-1}-1)(F_{n-2}-1)...(F_1-1),\quad\quad
n\geq 2.
\end{equation}
\end{cor}

\begin{ex}{\em
Let the sequence of finite cobweb posets $\{P_n\}_{n\geq 0}$ be
designated by the sequence $\{F_n\}_{n\geq 0}$ such that $F_n=n+1$
(i.e. by the sequence $\bf N$ of natural numbers). The examples of
corresponding characteristic polynomials are:\\\textrm{}\\
 $\chi_0(t)=1,$\\
 $\chi_1(t)=t-2,$\\
 $\chi_2(t)=t^2-2t+4,$\\
 $\chi_3(t)=t^3-2t^2+4t-18,$\\
 $\chi_4(t)=t^4-2t^3+4t^2-18t+120,$\\
 $\chi_5(t)=t^5-2t^4+4t^3-18t^2+120t-1050,$\\
 $\chi_5(t)=t^6-2t^5+4t^4-18t^3+120t^2-1050t+11340.$}
\end{ex}
\begin{ex}{\em
Let the sequence of finite cobweb posets $\{P_n\}_{n\geq 0}$ be
designated by the sequence $\{F_n\}_{n\geq 0}$ such that $F_1=1$
and $F_n=2n+1$ for $n\geq 1$. The examples of
corresponding characteristic polynomials are:\\\textrm{}\\
$\chi_0(t)=1,$\\
$\chi_1(t)=t-3,$\\
$\chi_2(t)=t^2-3t+10,$\\
$\chi_3(t)=t^3-3t^2+10t-56,$\\
$\chi_4(t)=t^4-3t^3+10t^2-56t+432,$\\
$\chi_5(t)=t^5-3t^4+10t^3-56t^2+432t-4224,$\\ }
\end{ex}
\begin{ex}{\em
Let the sequence of finite cobweb posets $\{P_n\}_{n\geq 0}$ be
designated by the sequence $\{F_n\}_{n\geq 0}$ such that $F_1=1$
and $F_n=k$ for $n\geq 1$ and for some $k>1$. The examples of
corresponding characteristic polynomials are:\\\textrm{}\\
$\chi_0(t)=1,$\\
$\chi_1(t)=t-k,$\\
$\chi_2(t)=t^2-kt+k(k-1),$\\
$\chi_3(t)=t^3-kt^2+k(k-1)t-k(k-1)^2,$\\
$\chi_4(t)=t^4-kt^3+k(k-1)t^2-k(k-1)^2t+k(k-1)^3.$\\

\noindent In general one has
$$\chi_n(t)=t^n-kt^{n-1}+k(k-1)t^{n-2}+...+(-1)^nk(k-1)^{n-1},\quad\quad n\geq 1.$$ }
\end{ex}

\noindent {\bf Acknowledgements}\\
Discussions with Participants of Gian-Carlo Rota Polish Seminar,\\
http://ii.uwb.edu.pl/akk/sem/sem\_rota.htm are highly appreciated.

 \end{document}